\documentclass{amsart}
\usepackage{graphicx}

\usepackage{amsfonts,amsmath,latexsym}
\usepackage{times}

\newtheorem{theorem}{Theorem}
\newtheorem{proposition}[theorem]{Proposition}

\def\Ss{\mathbb{S}}

\def\Rr{\mathbb{R}}
\def\dd{\mathrm{d}}

\def\ae{\langle}
\def\ad{\rangle}

\begin{document}

\title{Examples of scalar-flat hypersurfaces in $\mathbb{R}^{n+1}$}
\author{Jorge H. Lira \and Marc Soret}

\maketitle

\begin{abstract}
\emph{Given a  hypersurface $M$ of null scalar curvature in the unit
sphere $\mathbb{S}^n$, $n\ge 4$, such that its second fundamental
form has rank greater than 2, we construct a singular scalar-flat
hypersurface in $\Rr^{n+1}$  as a normal graph over a truncated cone
generated by $M$. Furthermore, this graph is 1-stable if the cone is
strictly 1-stable.}

\vspace{0.3cm}


\end{abstract}

\section{Introduction}

A consistent theme of research is the use of refined perturbation
techniques in the study of constant mean curvature surfaces and
metrics with positive constant scalar curvature. New and complex
examples and deep results on structure of moduli space of solutions
had been achieved with the aid of those techniques.

A kind of prototype of this type of construction may be found at the
seminal paper \cite{CHS}. There, the authors prove the existence of
minimal hypersurfaces with an isolated singularity in $\Rr^{n+1}$.
These examples arise as perturbations of cones over minimal
hypersurfaces of $\Ss^n$. 

Our contribution here focuses on a similar construction but for
scalar-flat singular hypersurfaces in Euclidean space $\Rr^{n+1}$.
We consider a truncated cone $\bar M^*$  in $\mathbb{R}^{n+1}$
generated by a hypersurface $M$ of $\mathbb{S}^n$ that satisfies
$S_2=0$ and then we take normal graphs over that cone. {\it A
priori} estimates plus a fixed point theorem assure the existence of
a graph with ``small'' boundary data which also satisfies the
equation $S_2=0$.

We recall that $S_2$ is one of the elementary symmetric functions
$S_r$, $1\le r \le n$, of the principal curvatures of a hypersurface
in
$\Rr^{n+1}$. 
An interesting feature of $S_2$ is that this curvature is {\it
intrinsic} and coincides with the scalar curvature of the
hypersurface.

Our aim here is to provide a test case that gives an evidence that
the well succeeded perturbation methods alluded above may be also
applicable to deal with some geometric problems involving fully
nonlinear elliptic equations. The results we obtained are in some
sense local. Global issues may be addressed only if we are able to
overcome serious technical difficulties.

\vspace{3mm} \noindent {\bf Theorem 1.} \emph{ Let $M$ be a
scalar-flat hypersurface in $\Ss^n$, $n\ge 4$. Suppose that the rank
of the second fundamental form of $M$ is greater than or equal to
$3$. Let $\psi$ be a function in $C^{2,\alpha}(M)$. There exists
$\Lambda<1$ depending on $M$ such that for each $\lambda\in
[0,\Lambda)$ there exists a function $u_{\lambda}$ defined in $\bar
M^*$ whose graph $\bar M^*_\lambda$  has null scalar curvature and
boundary given by $\Pi_J(u_\lambda)=\Pi_J(\lambda \psi)$, for some
integer $J$.}

\vspace{3mm}

\noindent Here, $\Pi_J$ is the projection map defined in p. 10.

This paper has the following presentation. In Section 2, we deduce
the null scalar curvature equation $\bar R(u)=0$  for the normal
graph of a function $u$ defined over $\bar M^*$. The linearized
equation involves the Jacobi operator $L$ in $\bar M^*$ which turns
to be elliptic in view of the hypothesis concerning the rank of the
second fundamental form of $M$. Section 3 is devoted to solve in
$\bar M^*$ a Dirichlet problem for the Jacobi operator with boundary
data $\psi$. Following closely \cite{CHS}, the idea is that an
adequate control of the data $f$ near the singular point in $\bar
M^*$ permits to solve $Lu=f$ in terms of separation of variables
techniques. Second order estimates for the resulting Fourier series
$u$ may be obtained in suitably weighted H\"older spaces. Applying
these estimates to the problem
\begin{equation}
\label{aux} Lu = Q(v),\quad u|_M = \psi,
\end{equation}
where $v$ is a function in a weighted H\"older space and $Q$
collects all nonlinear terms in $\bar R(v)=0$, we reduce the
nonlinear problem to that one of finding a fixed point for the map
that associates $v$ to the solution of (\ref{aux}). This is achieved
by showing that for small boundary data $\psi$, this map is a
contraction.

In the last section we relate the stability of the normal graphs
with the stability of the hypersurface $M\subset \mathbb{S}^n$.
There, stability refers to the functional $\mathcal{A}_1$ defined by
the integral of the mean curvature.

\vspace{3mm}

\noindent {\bf Theorem 2.} \emph{If $\bar M^*$ is strictly
$1$-stable, then the graph $\bar M^*_\lambda$ of the function
$u_{\lambda}$ given in Theorem 1 is strictly $1$-stable for
$\lambda$ sufficiently small.}

\vspace{3mm}

\noindent We point out that the results presented here may be easily
adapted to the other higher order mean curvatures $S_r, \, r\ge 3$.
It is interesting to produce examples with singular sets with small
codimension as Nathan Smale did for minimal hypersurfaces in
\cite{Smale}. This is the subject of current research by the
authors.

The corrections and suggestions by the anonymous referee improved
sensibly the reading of the paper. We express here our gratitude to
him.

\section{Scalar-flat cones}

\subsection{The scalar curvature equation.} Let $M$ be a compact hypersurface of the unit sphere $\Ss^{n}$ in
the Euclidean space $\Rr^{n+1}$. The {\it cone} over $M$ is the
hypersurface $\bar M$ in $\Rr^{n+1}$ parametrized by
\begin{equation}
X(t,\theta) =  t\,\theta, \quad t\in \Rr^+,\, \theta\in M.
\end{equation}
Let  $N$ be an unit normal vector field to $M$. Parallel
transporting $N$ along the rays $t\mapsto t\,\theta$ gives rise to a
normal vector field to $\bar M$. One then defines the first and
second fundamental forms of $\bar M$  respectively by
\begin{equation}
I = \ae \dd X,\dd X\ad, \quad II =-\ae \dd N,\dd X\ad.
\end{equation}
Let $x^1,\ldots, x^{n-1}$ be local coordinates in $M$ with
corresponding coordinate vector fields denoted by
$\partial_1,\ldots, \partial_{n-1}$. A local frame tangent to $\bar
M$ may be given by adding the vector field $\partial_t$ to that
coordinate  local frame. In terms of such a frame, the first
quadratic form is represented by the matrix
\begin{equation}
\big(\bar g_{\mu\nu}\big)=\left(\begin{array}{cc} t^2\,\theta_{ij} &
0
\\
0 & 1\end{array}\right)
\end{equation}
and the second fundamental form has components
\begin{equation}
\big(\bar b_{\mu\nu}\big)=\left(\begin{array}{cc} t\,b_{ij} & 0
\\
0 & 0\end{array}\right),
\end{equation}
where $\theta_{ij}=\ae \partial_i,\partial_j\ad$ and $b_{ij}=-\ae
\partial_j N,\partial_j\ad$ are the components of the first and second
fundamental forms of the immersion $M\subset\mathbb{S}^n$. Thus, the
Weingarten map $\bar A$ of $\bar M$ has local components given by
$\bar a^\mu_\nu = \bar g^{\mu\rho}\bar b_{\rho\nu}$. We then compute
\begin{equation}
\big(\bar a^\mu_\nu\big)=\left(\begin{array}{cc}
\frac{1}{t}\,a^i_{j} & 0
\\
0 & 0\end{array}\right),
\end{equation}
where $a^i_j = \theta^{ik}\, b_{jk}$ are the components of the
Weingarten map $A$ of $M$ defined by $N_i = - a^j_i \partial_j$.

If we denote by $\lambda_1,\ldots, \lambda_{n-1}$ the eigenvalues of
$A$, then the eigenvalues of $\bar A$ are
\begin{equation}
0,\frac{1}{t}\,\lambda_1,\ldots, \frac{1}{t}\,\lambda_{n-1}.
\end{equation}
The $r$-th mean curvature $\bar H_r$ of $\bar M$ is defined by
\begin{equation}
\bar H_r = \frac{1}{{n \choose r}}\,\bar  S_r,\quad 1\le r\le n,
\end{equation}
where $\bar S_r$ are the elementary symmetric functions of the
eigenvalues of $\bar A$ relative to $I$ given by
\begin{equation}
\det\big(\textrm{Id}-s\,\bar A\big) =1 - s\, \bar S_1 + s^{2}\,\bar
S_2+\ldots + (-s)^{n-1} \bar S_{n-1}+ (-s)^n\, \bar S_n.
\end{equation}
Denoting by $H_r$ and $S_r$  the corresponding functions on $M$, one
easily proves that
\begin{equation}
\bar S_r = \frac{1}{t^r}\, S_r, \quad 1\le r\le n-1
\end{equation}
and $\bar S_n =0$. For a given multi-index $i_1<\ldots <i_r$ with
$1\le i_k\le n$, we denote
\begin{equation}
D_{i_1\ldots i_r} = \det\big(\theta_{1j}\ldots b_{i_1j}\ldots
b_{i_rj}\ldots \theta_{n-1j}\big),
\end{equation}
that is, $D_{i_1\ldots i_r}$ is the determinant of the matrix
obtained  replacing in $(\theta_{ij})$ the columns numbered by
$i_1,\ldots,i_r$ by the corresponding columns in $(b_{ij})$.

In terms of these determinants, one calculates
\begin{equation}
\det(\theta_{ij})\, S_r =\sum_{i_1<\ldots<i_r}\, D_{i_1\ldots i_r}.
\end{equation}
We suppose that $M$ satisfies $S_2=0$. Thus, the cone $\bar M$ is a
scalar-flat manifold, that is, it holds that $\bar S_2=0$.

\subsection{The scalar curvature equation for normal graphs over cones.} From now on, we will be mainly concerned with
linearizing the equation $\bar S_2=0$ near $\bar M$. Given a
function $u:\bar M\to \Rr$ with sufficiently small $C^2$ norm, its
{\it normal graph}   is defined as the hypersurface
\begin{equation}
\bar M_u =\{ X(t,\theta) + u(t,\theta) \, N : t\in \Rr^+,\,\theta\in
M\}.
\end{equation}
We denote by $\bar S_2(u)$ the scalar curvature of $\bar M_u$. We
then proceed to linearize the equation $\bar S_2(u)=0$ and to
describe the nonlinear part of this equation.

We begin by determining the quadratic fundamental forms in $\bar
M_u$. The tangent space to $\bar M_u$ is spanned by the vector
fields $\theta + u_t \, N$ and
\begin{eqnarray}
t\,\big(\delta_i^j-u\,\bar a^j_i\big)\,\partial_j + u_i\,N,
\end{eqnarray}
where $u_t =\frac{\partial u}{\partial t}$ and $u_i =\frac{\partial
u}{\partial x^i}$. The induced metric in $\bar M_u$ has components
\begin{eqnarray*}
\bar g_{\mu\nu}(u)=\bar g_{\mu\nu} + \delta \bar
g_{\mu\nu},
\end{eqnarray*}
where
\begin{eqnarray*}
\big(\delta \bar g_{\mu\nu}\big)=\left(\begin{array}{cc} -2 u \bar
b_{ij}+ u^2 \bar r_{ij}
+ u_i u_j  & u_t u_i\\
u_i u_t &  u_t^2\end{array}\right)
\end{eqnarray*}
and $\bar r_{ij} =t^2\theta_{kl} \bar a_i^k  \bar a^l_j$ are the
components of the third fundamental form $\ae \dd N,\dd N\ad$ of
$\bar M$. More briefly, we may write
\begin{equation}
\delta\bar g_{\mu\nu} =- 2u\bar b_{\mu\nu}+u^2\bar r_{\mu\nu} +
u_\mu u_\nu.
\end{equation}
Let $\bar R_{\mu\nu}$ be the Ricci tensor of $\bar M$. If we denote
$\bar R = S_2$ and  $\bar R(u)= \bar S_2(u)$ then it follows that
\begin{eqnarray*}
\bar R (u) = \bar R + \delta \bar R,
\end{eqnarray*}
where
\begin{eqnarray*}
\delta \bar R  =  \bar g^{\mu\nu}\delta \bar R_{\mu\nu} +\delta \bar
g^{\mu\nu}\, \bar R_{\mu\nu}.
\end{eqnarray*}
A classical tensorial identity (see \cite{Bleecker}, p. 398) states
that
\begin{equation}
\label{deltag} \bar g^{\mu\nu}\,\delta \bar R_{\mu\nu}
=\bar\nabla_\rho W^\rho
\end{equation}
where $\bar\nabla$ denotes the Riemannian covariant derivative in
$\bar M$ with respect  to the metric $(\bar g_{\mu\nu})$ and
\begin{equation}
W^\rho =\bar g^{\rho\sigma}\bar g^{\mu\nu}\bar\nabla_\nu \delta \bar
g_{\mu\sigma}- \bar g^{\rho\nu}\bar g^{\mu\sigma}\bar\nabla_\nu
\delta \bar g_{\mu\sigma}.
\end{equation}
In what follows, we use the abbreviated notation $\bar\nabla^\rho
=\bar g^{\rho\mu}\bar\nabla_\mu$.

Since $\bar\nabla\bar g =0$ we may commute the covariant derivatives
and the components $\bar g^{\mu\nu}$ in the formula above
(\ref{deltag}), obtaining
\begin{eqnarray*}
\bar g^{\mu\nu}\,\delta \bar R_{\mu\nu} & = & \bar\nabla_\rho\bar
g^{\rho\sigma}\bar g^{\mu\nu}\bar\nabla_\nu \delta \bar
g_{\mu\sigma}- \bar\nabla_\rho \bar g^{\rho\nu}\bar
g^{\mu\sigma}\bar\nabla_\nu \delta \bar
g_{\mu\sigma}\\
&=&\bar\nabla_\rho\bar\nabla^\mu \bar g^{\rho\sigma}\delta \bar
g_{\mu\sigma}-\bar\nabla^\nu
\bar\nabla_\nu \bar g^{\mu\sigma}\delta \bar g_{\mu\sigma}\\
& = & -2\bar\nabla_\rho\bar\nabla^\mu \bar g^{\rho\sigma} \bar
b_{\mu\sigma}u+2\bar\nabla^\nu \bar\nabla_\nu \bar g^{\mu\sigma}
\bar b_{\mu\sigma}u +Q_1\\
& = & -2\bar\nabla_\rho\bar\nabla^\mu \bar
a^{\rho}_{\mu}u+2\bar\nabla^\nu \bar\nabla_\nu \bar a^{\mu}_{\mu}u
+Q_1\\
& =&-2\bar\nabla_\rho\bar\nabla^\mu \bar
a^{\rho}_{\mu}u+2\bar\nabla^\nu \bar\nabla_\nu \bar S_1\,u
+Q_1 \\
& = & 2\bar\nabla_\rho\bar\nabla^\mu \big(\delta^{\rho}_{\mu}\bar
S_1-\bar a^{\rho}_{\mu}\big)u+Q_1 \\
&=& 2\bar\nabla^\rho \bar\nabla_\mu \bar T^\mu_\rho \, u+Q_1,
\end{eqnarray*}
where $\bar T^\rho_\mu$ are the components of the $(1,1)$ tensor
field
\begin{equation}
\bar T_1 = \bar S_1 \,\textrm{Id}-\bar A
\end{equation}
and
\begin{eqnarray}
Q_1 &=& \bar\nabla_\rho\bar\nabla^\mu \big(u^2 \bar r^\rho_\mu +
u^\rho u_\mu\big)-\bar\nabla^\rho\bar\nabla_\rho\big(u^2 \bar
r^\mu_\mu + u^\mu u_\mu\big).
\end{eqnarray}
However, we have
\begin{eqnarray*}
& & \bar\nabla_\rho\bar\nabla^\mu \big(u^\rho u_\mu\big)
-\bar\nabla^\rho\bar\nabla_\rho\big(u^\mu u_\mu\big)=\bar
g^{\mu\nu}\bar g^{\rho\tau}\bar\nabla_{\rho}\bar\nabla_\nu(u_\tau
u_\mu)-\bar g^{\rho\tau}\bar g^{\mu\nu}
\bar\nabla_\tau\bar\nabla_\rho(u_\mu u_\nu)\\
& & \,\,\,\,=\bar g^{\mu\nu}\bar g^{\rho\tau}(u_{\tau;\nu\rho}
u_\mu+u_{\tau;\nu} u_{\mu;\rho}+u_{\tau;\rho} u_{\mu;\nu}+u_\tau
u_{\mu;\nu\rho})\\
& & \,\,\,\,\,\,\,\,\,-\bar g^{\rho\tau}\bar g^{\mu\nu}
(u_{\mu;\rho\tau} u_\nu+u_{\mu;\rho} u_{\nu;\tau}+u_{\mu;\tau}
u_{\nu;\rho}+u_\mu u_{\nu;\rho\tau})\\
& & \,\,\,\,=\bar g^{\mu\nu}\bar g^{\rho\tau}(u_{\tau;\nu\rho}
u_\mu+u_\tau u_{\mu;\nu\rho})-\bar g^{\rho\tau}\bar g^{\mu\nu}
(u_{\mu;\rho\tau} u_\nu+u_\mu u_{\nu;\rho\tau})\\
& & \,\,\,\,\,\,\,\,\,+\bar g^{\mu\nu}\bar g^{\rho\tau}(u_{\tau;\nu}
u_{\mu;\rho}+u_{\tau;\rho} u_{\mu;\nu})-\bar g^{\rho\tau}\bar
g^{\mu\nu} (u_{\mu;\rho} u_{\nu;\tau}+u_{\mu;\tau} u_{\nu;\rho}).
\end{eqnarray*}
Using Ricci identity
\[
u_{\nu;\tau\rho}- u_{\nu;\rho\tau}=\bar
R_{\tau\rho\sigma\nu}u^\sigma
\]
where $\bar R_{\tau\rho\sigma\nu}$ is the Riemann curvature tensor
in $\bar M$, we rewrite the terms with third order derivatives as
follows
\begin{eqnarray*}
& & \bar g^{\mu\nu}\bar g^{\rho\tau}(u_{\tau;\nu\rho} u_\mu+u_\tau
u_{\mu;\nu\rho})-\bar g^{\rho\tau}\bar g^{\mu\nu} (u_{\mu;\rho\tau}
u_\nu+u_\mu u_{\nu;\rho\tau})\\
& & \,\,\,\, =\bar g^{\mu\nu}\bar g^{\rho\tau}(u_{\tau;\nu\rho}
u_\mu-u_\mu u_{\nu;\rho\tau})+\bar g^{\mu\nu}\bar g^{\rho\tau}u_\tau
u_{\mu;\nu\rho}-\bar g^{\rho\tau}\bar g^{\mu\nu} u_{\mu;\rho\tau}
u_\nu
\\
& & \,\,\,\, =\bar g^{\mu\nu}\bar g^{\rho\tau}(u_{\nu;\tau\rho}
u_\mu- u_{\nu;\rho\tau}u_\mu)+\bar g^{\mu\nu}\bar g^{\rho\tau}
u_{\mu;\nu\rho}u_\tau-\bar g^{\rho\tau}\bar g^{\mu\nu}
u_{\mu;\rho\tau} u_\nu\\
& & \,\,\,\,=\bar g^{\mu\nu}\bar g^{\rho\tau}\bar
R_{\tau\rho\sigma\nu}u^\sigma u_\mu+\bar g^{\mu\nu}\bar g^{\rho\tau}
u_{\mu;\rho\nu}u_\tau+\bar g^{\mu\nu}\bar g^{\rho\tau} \bar
R_{\nu\rho\sigma\mu}u^\sigma u_\tau-\bar g^{\rho\tau}\bar g^{\mu\nu}
u_{\mu;\rho\tau} u_\nu\\
& & \,\,\,\,=\bar g^{\rho\tau}\bar R_{\tau\rho\sigma\nu} u^\sigma
u^\nu+\bar g^{\mu\nu}\bar g^{\rho\tau} u_{\mu;\rho\nu}u_\tau+\bar
g^{\mu\nu} \bar R_{\nu\rho\sigma\mu}u^\rho u^\sigma -\bar
g^{\rho\tau}\bar g^{\mu\nu} u_{\mu;\rho\tau} u_\nu.
\end{eqnarray*}
The antisymmetry of the curvature tensor in the last two indices
implies that $\bar g^{\rho\tau}\bar R_{\tau\rho\sigma\nu} u^\sigma
u^\nu=0$. Therefore, one has
\begin{eqnarray*}
& & \bar g^{\mu\nu}\bar g^{\rho\tau}(u_{\tau;\nu\rho} u_\mu+u_\tau
u_{\mu;\nu\rho})-\bar g^{\rho\tau}\bar g^{\mu\nu} (u_{\mu;\rho\tau}
u_\nu+u_\mu u_{\nu;\rho\tau})\\
& & \,\,\,\,=\bar g^{\mu\nu} u^{\rho} u_{\mu;\rho\nu}-\bar
g^{\rho\tau}u^{\mu} u_{\mu;\rho\tau} +\bar R_{\rho\sigma}u^\rho
u^\sigma\\
& & \,\,\,\,=\bar g^{\mu\nu} u^{\rho} u_{\mu;\rho\nu}-\bar
g^{\rho\tau}u^{\mu} u_{\rho;\mu\tau} +\bar R_{\rho\sigma}u^\rho
u^\sigma\\
& & \,\,\,\, =\bar R_{\rho\sigma}u^\rho u^\sigma.
\end{eqnarray*}
Thus, one concludes that
\begin{eqnarray*}
& & \bar\nabla_\rho\bar\nabla^\mu \big(u^\rho u_\mu\big)
-\bar\nabla^\rho\bar\nabla_\rho\big(u^\mu u_\mu\big) =\bar
R_{\rho\sigma}u^\rho u^\sigma+\bar g^{\mu\nu}\bar
g^{\rho\tau}(u_{\tau;\nu} u_{\mu;\rho}+u_{\tau;\rho}
u_{\mu;\nu})\\
& & \,\,\,\,\,\,\,\,-\bar g^{\rho\tau}\bar g^{\mu\nu} (u_{\mu;\rho}
u_{\nu;\tau}+u_{\mu;\tau} u_{\nu;\rho})\\
& & \,\,\,\,=\bar R_{\rho\sigma}u^\rho u^\sigma+u^\rho_{\,\,;\nu}
u^\nu_{\,\,;\rho}+u^\rho_{\,\,;\rho} u^\nu_{\,\,;\nu}-\bar
g^{\rho\tau} (u^\nu_{\,\,;\rho} u_{\tau;\nu}+u^\nu_{\,\,;\tau}
u_{\rho;\nu})\\
& & \,\,\,\, =\bar R_{\rho\sigma}u^\rho u^\sigma+u^\rho_{\,\,;\nu}
u^\nu_{\,\,;\rho}+u^\rho_{\,\,;\rho}
u^\nu_{\,\,;\nu}-u^\nu_{\,\,;\rho}
u^\rho_{\,\,;\nu}-u^\nu_{\,\,;\tau} u^\tau_{\,\,;\nu}\\
& & \,\,\,\,=\bar R_{\rho\sigma}u^\rho u^\sigma+u^\rho_{\,\,;\rho}
u^\nu_{\,\,;\nu}-u^\nu_{\,\,;\rho} u^\rho_{\,\,;\nu}.
\end{eqnarray*}
These calculations imply that
\begin{eqnarray*}
Q_1 & = & \bar\nabla_\rho\bar\nabla^\mu \big(u^2 \bar r^\rho_\mu
\big)-\bar\nabla^\rho\bar\nabla_\rho\big(u^2 \bar
r^\mu_\mu\big)+\bar R_{\rho\sigma}u^\rho u^\sigma+u^\rho_{\,\,;\rho}
u^\nu_{\,\,;\nu}-u^\nu_{\,\,;\rho} u^\rho_{\,\,;\nu}\\
& = & u^\rho_{\,\,;\rho} u^\nu_{\,\,;\nu}-u^\nu_{\,\,;\rho}
u^\rho_{\,\,;\nu} +2uu^\mu_{\,\,;\rho}(\bar
r^\rho_\mu-\delta^\rho_\mu \bar r)+u^\rho u^\mu (2\bar
r_{\rho\mu}-2\bar g_{\rho\mu}\bar r+\bar R_{\rho\mu})\\
& & \,\,\,\, +4uu^\rho(\bar r^\mu_{\rho;\mu}-\bar
r^\mu_{\mu;\rho})+u^2(\bar g^{\mu\nu}\bar r^\rho_{\mu;\nu\rho}-\bar
g^{\mu\nu}\bar r_{;\mu\nu}),
\end{eqnarray*}
where $\bar r = \bar r^\mu_\mu$.

It is a well-known fact that the tensor $\bar T_1$ is
divergence-free. Indeed, one computes using Codazzi's equation
\begin{equation*}
(\delta^\rho_\mu \bar S_1-\bar a^\rho_\mu)_{;\rho}= \delta^\rho_\mu
\bar a^\nu_{\nu;\rho}- \bar a^\rho_{\mu;\rho}= \bar a^\nu_{\nu;\mu}-
\bar a^\rho_{\rho;\mu}=0.
\end{equation*}
Using this, one gets
\begin{eqnarray*}
\bar g^{\mu\nu}\,\delta \bar R_{\mu\nu} &  = & 2\bar\nabla^\rho\big(
(\bar\nabla_\mu \bar T^\mu_\rho \big) u+\bar T^\mu_\rho
(\bar\nabla_\mu u)\big)+Q_1= 2\bar\nabla^\rho \bar T^\mu_\rho
\bar\nabla_\mu u +Q_1\\
&=& 2\bar\nabla_\rho \big(\bar T_\mu^\rho \bar\nabla^\mu u \big)+Q_1
=  2\textrm{div}\,\bar T_1\bar\nabla u+Q_1.
\end{eqnarray*}
On the other hand, we infer from Gauss equation that
\begin{eqnarray*}
\bar R_{\mu\nu}  =  \bar g^{\rho\sigma}\,\bar
R_{\mu\rho\nu\sigma}=\bar g^{\rho\sigma}\big(\bar b_{\mu\nu}\bar
b_{\rho\sigma}-\bar b_{\mu\sigma}\bar b_{\nu\rho}\big)=\bar
b_{\mu\nu}\,\bar S_1-\bar r_{\mu\nu}
\end{eqnarray*}
and since
$$\delta \bar g^{\mu\nu}\bar R_{\mu\nu}= \delta\bar
g^{\mu\nu}g_{\mu\rho}\bar R^\rho_\nu=- \delta g_{\mu\rho}\bar
g^{\mu\nu}\bar R^\rho_\nu= - \delta g_{\mu\rho}\bar R^{\mu\rho},
$$
one obtains
\begin{eqnarray*}
\delta \bar g^{\mu\nu}\, \bar R_{\mu\nu} & = &  2u \bar S_1 \bar
b_{\mu\nu}\bar b^{\mu\nu} -2u\bar b_{\mu\nu}\bar
r^{\mu\nu}+Q_2=2u\bar S_1 \textrm{tr} \bar A^2 -2u\textrm{tr} \bar
A^3 +Q_2\\
&=& 2\textrm{tr}\big((\bar S_1 \textrm{Id}-\bar A)\bar A^2\big)+Q_2=
2\textrm{tr}\bar T_1\bar A^2+Q_2,
\end{eqnarray*}
where
\begin{equation}
Q_2 = -\bar R^{\mu\rho}\big(u^2\bar r_{\mu\rho}+ u_\rho u_\mu\big).
\end{equation}
Since we are assuming that $\bar S_2=0$ one easily verifies that
\begin{equation}
\textrm{tr}\,\bar T_1\bar A^2 = -3\bar S_3.
\end{equation}
We then conclude that the equation $\bar R(u)=0$ may be written as
\begin{equation}
\label{nonl-eqtn} Lu + Q(u) =0,
\end{equation}
where
\begin{equation}
L u =\textrm{div} \,\bar T_1\bar\nabla u-3\bar S_3 u
\end{equation}
is the Jacobi operator for the scalar curvature and  $Q=Q_1+Q_2$.


The quadratic term $Q$ has the form
\begin{eqnarray}
\label{Q} Q(u,\bar\nabla u,\bar\nabla^2 u) =|\Delta_{\bar M}
u|^2-|\bar\nabla^2 u|^2+\textrm{tr}\, P_0(u)\cdot\bar\nabla^2
u+P_1(u,\bar\nabla u)
\end{eqnarray}
where $\Delta_{\bar M}$ is the Laplace-Beltrami operator in $\bar M$
and
\begin{eqnarray*}
(P_0)^\rho_\mu=2(\bar r^\rho_\mu-\delta^\rho_\mu \bar r)u
\end{eqnarray*}
and
\begin{eqnarray*}
& & P_1 = (2\bar r_{\rho\mu}-2\bar g_{\rho\mu}\bar r)u^\rho
u^\mu+4(\bar r^\mu_{\rho;\mu}-\bar r^\mu_{\mu;\rho})uu^\rho\\
& & \,\,\,\,+(\bar g^{\mu\nu}\bar r^\rho_{\mu;\nu\rho}-\bar
g^{\mu\nu}\bar r_{;\mu\nu}-\bar R^{\mu\nu}\bar r_{\mu\nu})u^2.
\end{eqnarray*}


\section{The Dirichlet problem for the Jacobi operator.} As we proved above,
a normal graph $\bar M_u$ is scalar-flat if $u$ satisfies the fully
nonlinear equation (\ref{nonl-eqtn}). Our goal in this section is to
solve the corresponding linearized equation for small boundary data
by using Fourier analysis in some suitably weighted spaces.

Following the notation previously fixed, we denote
\begin{equation}
\bar L_1 u = \textrm{div} \,\bar T_1\bar\nabla u.
\end{equation}
The corresponding tensor and operator in $M$ are respectively
\[
T_1 = S_1 \textrm{Id}- A
\]
and
\begin{equation}
L_1 u = \textrm{div}\, T_1\nabla u,
\end{equation}
where the divergence and gradient are taken this time on $M$. In
\cite{BdC}, it is proved that the operators $L$ and $\bar L_1$
decomposes as follows
\begin{equation}
\bar L_1 u = \frac{1}{t}\, S_1  u_{tt} + \frac{n-2}{t^2} \, S_1 u_t
+ \frac{1}{t^3} \,L_1 u(t,\cdot)
\end{equation}
and
\begin{equation}
\label{Jtt} L u = \frac{1}{t}\, S_1  u_{tt} + \frac{n-2}{t^2} \, S_1
u_t + \frac{1}{t^3} \,L_1 u(t,\cdot)-3\frac{1}{t^3}S_3 u.
\end{equation}
From now on, we assume that $S_3$ never vanishes along $M$ or
equivalently that $\textrm{rk}A\ge 3$. In \cite{HL}, it is proved
that this assumption assures the ellipticity of the second-order
differential operator $L$. This is a crucial ingredient in our
analysis. We point out that there are examples of hypersurfaces
fitting our assumptions in $\Ss^n$ like certain products of spheres.

As an example, if we fix the lowest dimension $n=4$, we may consider
the product of spheres $M=\mathbb{S}^2(a_1)\times\mathbb{S}^1(a_2)$
immersed in $\mathbb{S}^4$, where $a_1 =\sqrt{1/3}$ and $a_2
=\sqrt{2/3}$. With these choices one has $S_2 =0$ and
\[
S_1 = 2\sqrt{2} -\sqrt{1/2}\quad \textrm{and}\quad S_3 =-\sqrt{2}.
\]
For a detailed explanation on these products of spheres, we refer
the reader to \cite{BdC}.

We begin our analysis of the equation (\ref{nonl-eqtn}) by solving
first the non-homogeneous linear Dirichlet problem for the Jacobi
operator
\begin{eqnarray}
\label{jlinear} 
L u   =  f \,\, \textrm{in}\,\,\bar M^*,\quad u = \psi \,\,
\textrm{in}\,\, M
\end{eqnarray}
where $\bar M^*$ is the truncated cone obtained restricting the
variable $t$ to $(0,1]$. Using (\ref{Jtt}), we reduce the linear
equation $Lu = f$  to
\begin{equation}
\label{radial} t^2 S_1 u_{tt} + (n-2) t S_1  u_t +L_1 u(t,\cdot) -
3S_3 u  = t^3 f(t, \cdot).
\end{equation}
The hypothesis on $S_3$ implies that  $S_1$ also never vanishes. We
then  may choose an orientation for $M$ in such a way that $S_1 >
0$. Hence, the operator in $M$ defined by
\begin{equation}
-S_1^{-1}(L_1 -3S_3)
\end{equation}
has $L^2(M,S_1\dd \theta)$ discrete spectra given by a set of
diverging eigenvalues
\begin{equation}
\mu_1\le \mu_2\le\ldots\to +\infty
\end{equation}
with corresponding eigenfunctions $\{\phi_1,\phi_2,\ldots\}$. These
facts permit to separate variables in (\ref{radial}) and reduce the
problem to the determination of a Fourier series for $u$.
We will see that a formal solution of (\ref{radial}) in Fourier
series gives rise to convergent solutions if we  consider functions
$f= f(t,\theta)$ such that
\begin{equation}
\label{hyp1}
 \vert f\vert _t ^2 : = \int_{M} f(t,\theta)^2
S_1(\theta)^{-1} \dd \theta <\infty,\quad t\in (0,1].
\end{equation}
Let $m>2$ and $\epsilon>0$ be real constants to be chosen later. It
is required too that the function $t\mapsto \vert f\vert_t$
satisfies 
\begin{equation}
\sup_{(0,1)} t^{2 -m -\epsilon}  \vert f\vert_t <\infty.
\end{equation}
This implies that $f(0,\cdot) =0$ and
\begin{equation}
\label{hyp2} ||f||: =\left( \int^1_0 t^{4-2m} \vert f\vert_t^2 \dd t
\right)^{1/2}<\infty.
\end{equation}
Under the  assumptions above on $f$, it is possible to decompose it
in its Fourier series
\begin{equation}
\frac{f}{S_1} = \sum_{j=1}^\infty f_j(t) \phi_j(\theta)
\end{equation}
with $f_j(t) = \int f\phi_j\,\dd \theta$. Let $u$ be a formal
solution
\begin{equation}
u(t,\theta) = \sum_j a_j(t) \phi_j(\theta)
\end{equation}
of equation (\ref{radial}). Thus, the coefficients $a_j$ are
determined by the sequence of ODE's
\begin{equation}
\label{ode} t^2 a_j'' + (n-2) t a_j' - \mu_j a_j = t^3 f_j,\quad j=
1,2,\ldots
\end{equation}
The  homogeneous equations associated to (\ref{ode}) have solutions
of the form $t^{\gamma_j}$ where  $\gamma_{j}$ is root of the
characteristic equation $\gamma^2 + (n-3)\gamma - \mu_j = 0$. Its
roots are the {\it indicial roots}
\begin{equation}
\label{iroots} \gamma_j = -\frac{n-3}{2} \pm \sqrt{
\frac{(n-3)^2}{4}+\mu_j}.
\end{equation}
We observe that $\gamma_j$ may be complex since $\mu_j$ may be
negative. In these cases, one has $\Re \gamma_j = (3-n)/2$. Since
the eigenvalues $\mu_j$ diverge to $+\infty$, there
exists 
an index $J$  such that   $\Re (\gamma_{J+1}) = \gamma_{J+1}
>0$. This index may be chosen so that for a given  $m>2$ it holds that
\begin{equation}
\frac{3-n}{2}\le\ldots\le \Re (\gamma_J) < m < \Re (\gamma_{J+1})\le
\Re (\gamma_{J+2})\le...
\end{equation}
From now on, we consider these choices for $m$ and $J$.

In order to find a particular solution of the non-homogeneous
equation (\ref{ode}), we consider functions of the form $a_j(t) =
t^{\gamma_j}v_j(t)$. Plugging this expression of $a_j$ in
(\ref{ode}) we obtain
\begin{equation}
t^{\gamma_j+2}v_j'' + (2\gamma_j + n - 2) t^{\gamma_j+1}v_j' = t^3
f_j
\end{equation}
and after multiplying this equation by $t^{\gamma_j+ n-4}$ one has
\begin{equation}
(t^{n-2 +2\gamma_j}v_j')' = t^{n-1+\gamma_j} f_j.
\end{equation}
Integrating twice we get
\begin{equation}
v_j = \alpha_j +\int^t_{\beta_j} s^{2-n-2\gamma_j}\int^s_0
\tau^{n-1+\gamma_j} f_j\,\dd\tau\,\dd s, \quad j=1,2,\ldots
\end{equation}
where $\alpha_j$ and $\beta_j$ are constants of integration to be
specified in the sequel. We conclude that the formal solution
$u=\sum_j a_j \,\phi_j$ to equation (\ref{radial})  has coefficients
of the form
\begin{equation}
\label{solution}
 a_j(t) = \Re\bigg(\alpha_j t^{\gamma_j}+
 t^{\gamma_j}\int^t_{\beta_j}s^{2-n-2\gamma_j}\,\int^s_0
 \tau^{n-1+\gamma_j} f_j(\tau)\,\dd\tau\,\dd s\bigg).
 \end{equation}
We claim that the integrals in the definition of these coefficients
converge in $(0,1]$ if we choose $\alpha_j = \beta_j =0$ for $j\le
J$ and $\beta_j=1$ for $j\ge J+1$. In fact, one has
\begin{eqnarray*}
f_j(t) = \int_M \frac{f}{\sqrt{S_1}}\phi_j \sqrt{S_1} \,\dd\theta
\le \sqrt{\int_M \frac{f^2}{S_1}\,\dd\theta} \,\,\sqrt{\int_M
\phi^2_j S_1\,\dd\theta}=|f|_t.
\end{eqnarray*}
Thus, using the hypothesis (\ref{hyp2}) and Cauchy-Schwarz
inequality we estimate, for a constant $c$ that does not depend on
$f$,
\begin{eqnarray*}
\int^s_0 \tau^{n-1+\gamma_j}f_j(\tau)\,\dd\tau &\le & \sqrt{\int^s_0
\tau^{2(n-1+\gamma_j)}\tau^{2m-4}\, \dd\tau}
\,\,\sqrt{\int^s_0 \tau^{4-2m}|f|_\tau^2\,\dd\tau}\\
& = & c||f||s^{n-3+m+\gamma_j+\frac{1}{2}},
\end{eqnarray*}
where we used the fact that $m>\Re{\gamma_j}$ for $j\le J$ in order
to assure convergence of the integral at $s=0$. This estimate
implies that
\begin{eqnarray}
\label{est-aj} t^{\gamma_j}\int^t_{\beta_j}
s^{2-n-2\gamma_j}\,\int^s_0 \tau^{n-1+\gamma_j}f_j(\tau)\,
\dd\tau\,\dd s \le c||f|| \,t^{\gamma_j}\int^t_{\beta_j}
s^{m-\gamma_j-\frac{1}{2}}\, \dd s.
\end{eqnarray}
For $j\le J$, the right hand side converges at $t=0$ if one sets
$\beta_j=0$. For $j\ge J+1$, it converges if we consider $\beta_j
=1$. This proves the claim.

The values of $\alpha_j$ for $j\ge J+1$ are determined by
\begin{eqnarray}
\label{boundary} \alpha_j :=\int_{M}\lim_{t\to 1}
u(t,\cdot)\,\phi_j,\quad j\ge J+1.
\end{eqnarray}
Let $\Pi_J$ be the projection of $L^2(M,S_1\dd\theta)$ in the linear
subspace spanned by the eigenfunction $\phi_j,\,j\ge J+1$. Thus,
\begin{equation}
\Pi_J (u) = \Pi_J(\psi)
\end{equation}
if and only if
\begin{equation}
\Pi_J(\psi) = \sum_{j=J+1}^\infty \alpha_j \phi_j.
\end{equation}
Thus, since $\psi\in L^2(M,S_1\dd\theta)$, one has
\begin{equation}
\sum_{j= J+1}^\infty \alpha_j^2 <\infty.
\end{equation}
In this case, we then had verified that the problem (\ref{jlinear})
has as solution the convergent Fourier series $u$ defined by the
coefficients $a_j$ above.

In particular we have found a solution to the equation $Lu=0$ with
boundary Dirichlet data $\psi$ referred to in what follows as the
$L$-harmonic extension of $\psi$. In other terms we denote by $H_J
(\psi)$ the Fourier series solution of
\[
Lu=0 \,\,\textrm{ in }\,\, \bar M^*, \quad \Pi_J(u) = \Pi_J
(\psi)\,\,\textrm{ in }\,\, M.
\]
Notice that our previous calculations imply that
\begin{equation}
\label{HJ} H_J(\psi) = \sum_{j=J+1}^\infty \alpha_j
t^{\gamma_j}\phi_j.
\end{equation}
and $H_J$ is a right inverse to $\Pi_J$.

In order to obtain integral estimates for $u$, we notice that since
\begin{equation}
|uS_1|_t^2 =\int_M u^2(t,\theta)S_1(\theta)\,\dd\theta =
\sum_{j=1}^\infty a_j^2(t)
\end{equation} it follows that
\begin{equation}
|u|_t^2\le c \sum_{j=1}^\infty a_j^2(t)
\end{equation}
where $c= 1/(\inf_M S_1^2(\theta))$. On the other hand, using
(\ref{solution}) and (\ref{est-aj}), one obtains from Cauchy-Schwarz
inequality
\begin{equation}
t^{-m}|u|_t \le t^{-m}\sqrt{\sum_{j=1}^\infty a_j^2(t)}\le c||f||
+\sqrt{\sum_{j=J+1}^\infty \alpha_j^2},
\end{equation}
where $c>0$ is a positive constant which depends on $M, m$ and $J$.
In a similar way, using (\ref{solution}) and (\ref{HJ}) one proves
that
\[ t^{-m}|u-H_J(\psi)|_t \le c||f||.
\]

We summarize the facts above in the following proposition.

\vspace{0.3cm}

\begin{proposition}
\label{C0} Let  $m>2$ be a constant and let  $J$ be an integer such
that
$$0< \Re(\gamma_J) < m < \Re(\gamma_{J+1})$$
for $\gamma_j$ given by (\ref{iroots}). Given a function $f$ defined
in $\bar M^*$ satisfying
\[
\sup_{(0,1)}t^{2-m-\epsilon}|f|_t <\infty
\]
and a function $\psi\in L^2(M,S_1\dd\theta)$, the series
$$u= \sum_{j=1}^\infty  a_j\, \phi_j$$
with $a_j$ defined by (\ref{solution}) is the unique solution of
\begin{equation}
\label{jpilinear} Lu= f \textrm{ in }
 \bar M^*\quad \textrm{ and }\quad
\Pi_J(u) = \Pi_J(\psi) \textrm{ in } M
\end{equation}
satisfying
\begin{equation}
\label{ut} \sup_{(0,1)} t^{-m}|u|_t<\infty.
\end{equation}
Moreover, we have the following estimates for $u$
\begin{eqnarray}
& & t^{-m}|u|_t \le c\,(||f|| + |\Pi_J(\psi)|),\nonumber\\
& & t^{-m}|u-H_J(\psi)|_t \le c\,||f||,
\end{eqnarray}
where the constant $c$ does not depend on  $f$.
\end{proposition}

\vspace{0.3cm}




\noindent{\it Proof of the uniqueness.} In view of the previous
discussion, it remains to prove the uniqueness of the solution. If
we consider two solutions $u_1$ and $u_2$ of the equation $Lu =f$,
then their difference $v= u_1 - u_2$ is decomposed as $v = \sum_j
b_j \,\phi_j$ where the functions $b_j$ are solutions of the
homogeneous ODE associated to (\ref{ode}). Notice that (\ref{ut})
implies that
 $u_1$ and $u_2$ vanish at the origin. Thus,
$b_j(t)\to 0$ as $t\to 0$ for all $j$. Moreover, if $j\ge J+1$ then
$\gamma_j$ is real and positive. So, $\mu_j$ is necessarily
positive. Therefore the maximum principle guarantees that $b_j = 0$
for all $j\ge J+1$. For $j\le J$ we have that $b_j$ is of the form
$b_j = c t^{\gamma_j} + \tilde c t^{\tilde \gamma_j}$ where
$\gamma_j, \, \tilde \gamma_j$ are the roots of the characteristic
equation. Thus $|t^{-m}b_j |\to \infty$ unless that $b_j = 0$ for
$j\le J$. So, we have proved the proposition.

\vspace{0.3cm}


\noindent Following \cite{CHS} we now define some weighted H\"older
spaces in terms of that it is possible to obtain second order
estimates for the solution of the linear problem.

More precisely, we introduce as in \cite{CHS} and \cite{PR}, the
norms
\begin{eqnarray}
\label{pnorms} & & |v|_{k,\alpha,t}= \sum_{l=0}^k t^l |\bar\nabla^l
v |_{0,A_t}+t^{k+\alpha}[\bar\nabla^k u]_{\alpha, A_t},
\end{eqnarray}
for $t\in (0,1/2)$, $k$ a positive integer and $\alpha\in (0,1)$.
Here, $A_t$ is the truncated cone corresponding to $t< |X| < 2t$ and
$|\cdot|_{0,\alpha,A_t}$ denotes the usual H\"older norm in $A_t$.

\vspace{0.3cm}
\begin{proposition}
\label{C2} Under the hypothesis of the Proposition \ref{C0}, the
function $u$ satisfies
\begin{eqnarray}
\label{C2est} & & t^{-m}|u|_{2,\alpha,t}\le c\,(||f||_{\alpha} +
|\Pi_J(\psi)|),\nonumber\\
& & t^{-m}|u - H_J(\psi)|_{2,\alpha,t} \le c\, ||f||_{\alpha},
\end{eqnarray}
for $t\in (0,\frac{1}{2})$, $\psi\in C^{2,\alpha}(M)$ and
\begin{equation}
||f||_{\alpha}\equiv \sup_{0<t<1/2}t^{2-m-\epsilon}|f|_{0,\alpha,t}
\end{equation}
where $\epsilon$ is a fixed positive number. The constants do not
depend on $f$.
\end{proposition}
\vspace{0.3cm}

\noindent {\it Sketch of the proof.} A similar estimate for the
Laplacian could be found in \cite{P}  and \cite{PR}. We may obtain
the estimates for elliptic linear operators with constant
coefficients and  only second order terms. The general case could be
handled by freezing coefficients in $L$. For usual H\"older norms,
this method is nicely exposed   in Chapters 4 and 6 of  \cite{GT}.

\vspace{0.3cm}

\section{Solving the nonlinear problem}



Using the weighted H\"older spaces we just defined above, we then
introduce the subspace $B$ of $C^{2,\alpha}(\bar M^*)$ consisting of
the functions $v$ for which
\begin{equation}
\label{norm} ||v|| = \sup_{0<t<1/2} t^{-m}|v|_{2,\alpha,t}
\end{equation}
is finite.

We define a map $U$ in the unit ball in $B$ in the following way:
given a function $v\in B$ with $||v||<1$, $U(v)$ is the solution of
the linear problem
\[
L U = Q(v) \,\, \textrm{ in  }\,\, M^*, \quad \Pi_J(U) = \Pi_J
(\psi)\,\,\textrm{ in }\,\, M
\]
as defined in Proposition 1. Our task now is to exhibit a convex
subset $K$ of the unit ball in $B$ so that $U|_K$ is a contraction
map.

With this purpose, we begin by estimating $Q(v)$ for  $v$ with
$||v||<1$. We have, using that $t<1$,
\begin{eqnarray*}
& & |Q(v)|_{0,\alpha,t}\le 2|\bar\nabla^2 v |^2_{0,\alpha,t} +
|P_0|_{0,\alpha,t}|\bar\nabla^2 v|_{0,\alpha,t} +
|P_1|_{0,\alpha,t} \\
& & \,\,\,\,\le 2 (t^{-2}|v|_{2,\alpha,t})^2 +
|P_0|_{0,\alpha,t}t^{-2}|v|_{2,\alpha,t} + |P_1|_{0,\alpha,t}\\
& & \,\,\,\,\le  2 t^{-4}|v|_{2,\alpha,t}^2 +
C_0t^{-2}|v|_{0,\alpha,t}|v|_{2,\alpha,t} + C_1 (|v|_{0,\alpha,t} + |\bar\nabla v|_{0,\alpha,t})^2\\
& & \,\,\,\,\le 2  t^{-4}|v|_{2,\alpha,t}^{2} + C_0 t^{-2}
|v|^2_{2,\alpha,t}+ C_1
(1+t^{-1})^2|v|^2_{2,\alpha,t} \\
& & \,\,\,\,\le \mu |v|_{2,\alpha,t}^2 \le \mu t^{2m}||v||^2,
\end{eqnarray*}
where $C_0, C_1$ and $\mu$ are positive constants depending only on
$M$.

We choose $\epsilon$ such that $m+2\ge \epsilon$. Since $t<1$ we
have $t^{2m}\le t^{m-2+\epsilon}$. Thus we obtain
\begin{equation}
\label{Bnorm} |Q(v)|_{0,\alpha,t}\le \mu t^{m-2+\epsilon}||v||^2
\end{equation}
and similarly one easily verifies that
\begin{equation}
\label{BHnorm} |Q(v) - Q(w)|_{0,\alpha,t}\le \mu (||v||+ ||w||)(||v
- w||)t^{m-2+\epsilon}.
\end{equation}
It follows from estimates stated in   Proposition \ref{C2} that
$U(v)$
satisfies
\begin{eqnarray*}
||U(v) - H_J\psi|| &=& \sup_{0<t<1/2} t^{-m}|U(v) -
H_J\psi|_{2,\alpha,t}\le c
||f||_{\alpha}\\
&  = & c \sup_{0<t<1/2}t^{2-m-\epsilon}|Q(v)|_{0,\alpha,t}\le c\mu
||v||^2.
\end{eqnarray*}
Moreover since  $L(U(v)-U(w))= Q(v) -Q(w)$ and $\Pi_J (U(v))=
\Pi_J(U(w))$ then using the first estimate in Proposition 2 we
obtain
\begin{eqnarray*}
||U(v) -U(w)|| & = &  \sup_{t} t^{-m}|U(v) -U(w)|_{2,\alpha,t}\le
c||Q(v)-Q(w)||_{\alpha}\\ & = & c \sup_{t}
t^{2-m-\epsilon}|Q(v)-Q(w)|_{0,\alpha,t} \\
& \le & c\mu (||v|| +||w||)(||v -w||).
\end{eqnarray*}
In view of the last inequality, it is necessary to distinguish two
cases. We suppose first  that $c\mu <\lambda/2$ for some constant
$\lambda<1$. Then, given $u,v$ with $||u||\le 1$ and $||v||\le 1$ we
have
$$||U(u) - U(v)|| \le \lambda
||u-v||.$$ Moreover,
$$||U(v)|| \le c\mu ||v||^2 + ||H_J\psi||\le 1$$
if we assume that
$$||H_J\psi||\le 1-c\mu ||v||^2.$$
Since $||v||\le 1$ the last inequality holds if we suppose
\begin{equation}
\label{psi-1} ||H_J\psi||\le 1-c\mu,
\end{equation}
which is true for suficiently small $\psi$. Hence, assuming this we
conclude that $U|_K : K\to K$ is a contraction map where $K$ is the
intersection of the unit open ball in $B$ with the affine subspace
$\mathcal{P}=\{v\in B:\Pi_J v = \Pi_J \psi\}$. Notice that the
smallness of $\psi$ also guarantees that $K$ is not empty.

Now, we suppose that $c\mu \ge 1/2$. In this case,  we assume that
$||v||\le a$ for some constant $a$ to determine. One gets
$$||U(v)|| \le c\mu ||v||^2 +||H_J\psi|| \le c\mu a^2 +
||H_J\psi||.$$ Thus in order that $||U(v)|| \le a$ it is sufficient
that
$$c\mu a^2 -a + ||H_J\psi||\le 0.$$
Then $a$ must be choosen as $a\le \frac{1+\sqrt{1-4c\mu
||H_J\psi||}}{2c\mu}$. We must assume that
$$||H_J\psi||\le \frac{1}{4c\mu}$$
in order to assure that the square root above is well-defined. Since
$$\frac{1}{2c\mu}< \frac{1+\sqrt{1-4c\mu
||H_J\psi||}}{2c\mu},$$ we may choose $a=1/(2c\mu)$. So, we must
suppose simultaneously that $||v||\le 1$ and that $||v||\le a$.
However, the hypothesis $c\mu \ge 1/2$ implies that $a=1/(2c\mu) \le
1$. So, we prove that $U(K_1) \subset K_1$ and $U|_{K_1}$ is a
contraction mapping, where $K_1$ is the intersection of the ball of
radius $a$ in $B$ with the affine plane $\mathcal{P}$.

In both cases,  we had just  verified  that $U$ defines a
contraction map in properly chosen convex sets of the Banach space
$B$. So, by  Leray's fixed point theorem (see, e.g., \cite{GT},
Chapter 11), we assure the existence of a solution for the equation
(\ref{nonl-eqtn}).

\begin{theorem} Let $M$ be a scalar-flat hypersurface in $\Ss^n$, $n\ge 4$. Suppose that the
rank of the second fundamental form of $M$ is greater than or equal
to $3$. Let $\psi$ be a function in $C^{2,\alpha}(M)$. There exists
$\Lambda<1$ depending on $M$ such that for each $\lambda\in
[0,\Lambda)$ there exists a function $u_{\lambda}$ defined in $\bar
M^*$ such that the graph $\bar M^*_\lambda$ of $u_{\lambda}$ has
null scalar curvature and boundary given by
$\Pi_J(u_\lambda)=\Pi_J(\lambda \psi)$, for some integer $J$.
\end{theorem}


\section{Stability of scalar-flat cones}

It is well-known that scalar-flat hypersurfaces in $\Rr^{n+1}$ are
locally characterized as extrema of the action
\begin{equation}
\mathcal{A}_1 = \int_{\bar M} \bar S_1\, \dd \bar M.
\end{equation}
In this context, the Jacobi operator $L$ is naturally linked to
stability of the hypersurface. For details, we refer the reader to
\cite{Re}, \cite{Ro} and \cite{BC}.

In this section, we are concerned with the stability of the
scalar-flat cones and graphs we had defined above. For that, we
consider a function $u\in C^2_0(\bar M^*)$. The first and second
variation formulae for $\mathcal{A}_1$ are:
\begin{eqnarray*}
& & \mathcal{A}_1'(0)= 0,\quad \mathcal{A}_1 ''(0) = -\int_{\bar
M^*} u\, Lu \,\mathrm{d}\bar M.
\end{eqnarray*}
We recall that the Jacobi operator in the last formula is
\begin{eqnarray*}
Lu =\bar L_1u - 3\bar S_3 u = S_1 t^{1-n} (t^{n-2}u_t)_t+
\frac{1}{t^3}\,(L_1 u(t,\cdot)-3S_3u).
\end{eqnarray*}
We decompose $u$ in its Fourier coefficients with respect to the
eigenfunctions $\{\phi_j\}$ of $-\frac{1}{S_1}\big(L_1-3S_3)$
obtaining $u = \sum_j b_j \phi_j$ with $b_j(0)=b_j(1)=0$ and
\begin{eqnarray*}
& & L u  = \sum_j S_1\,\big( t^{1-n}(t^{n-2}b_j')_t-t^{-3}\mu_j
b_j\big)\,\phi_j.
\end{eqnarray*}
Since the metric of $\bar M^*$ in spherical coordinates $(t,\theta)$
is written in the form $\mathrm{d}t^2 + t^2
\theta_{ij}\mathrm{d}\theta^i\otimes \mathrm{d}\theta^j$, one has
$\dd \bar M  = t^{n-1}\dd t \,\dd \theta$, where $\dd \theta$ is the
volume form in $M$. Since $b_j(1)=0$, for all $j$, it results that
\begin{eqnarray*}
 \int_{\bar M^*} u \, L u \,\dd \bar M & = & \sum_{j,k}\int^1_0
((t^{n-2}b_j')_t-t^{n-4}\mu_j b_j)b_k \int_{M} \phi_j\,\phi_k
S_1(\theta) \mathrm{d}\theta
\\
&  = & \int^1_0 \sum_j((t^{n-2}b_j')_tb_j-t^{n-4}\mu_j
b_j^2) \,\mathrm{d}t\\
&   = &  -\int^1_0 \sum_j (t^{n-2}(b_j')^2 + t^{n-4}\mu_j
b_j^2)\,\mathrm{d}t.
\end{eqnarray*}
The first term in the last integral is given by
\begin{equation}
\int_{\bar M^*} u_t^2 \bar S_1 \dd \bar M= \int_{\bar M^*} u_t^2
t^{-1}S_1\dd \bar M= \int^1_0 t^{n-2}\sum_j (b_j')^2 \dd t.
\end{equation}
Denote $\mu_1^-=\max\{-\mu_1,0\}$, where $\mu_1$ is the smallest
eigenvalue of the operator $-\frac{1}{S_1}(L_1-3S_3)$. Thus, one
obtains
\begin{eqnarray*}
\label{ineq-1} & & \int_{\bar M^*} u \,L u\,\dd \bar M \le
-\int_{\bar M^*} u_t^2 \bar S_1\dd \bar M + \mu^{-}_1 \int^1_0
t^{n-4} \sum_ j b_j^2\,\mathrm{d}t.
\end{eqnarray*}
However, one has
\begin{equation}
\int_{\bar M^*} u^2 t^{-2}\bar S_1 \dd \bar M= \int_{\bar M^*} u^2
t^{-3}S_1 \,\dd\bar M = \int^1_0 t^{n-4} \sum_j b_j^2 \,\dd t
\end{equation}
and the expression on the right hand side of (\ref{ineq-1}) may be
calculated as follows
\begin{eqnarray*}
\int^1_0 t^{n-4}\sum_j b_j^2 \mathrm{d}t &=& \frac{1}{n-3}\int_0^1
(t^{n-3}\, \sum_j b_j^2)_t\,\mathrm{d}t
  -\frac{2}{n-3}\,\int^1_0 t^{n-3}
\sum_j b_j b_j' \dd t\\
&  \le & \frac{2}{n-3}\bigg( \int^1_0 t^{n-2}\sum_j (b_j')^2  \dd
t\bigg)^{1/2} \bigg(\int^1_0 t^{n-4}\sum_j b_j^2 \dd t\bigg)^{1/2}.
\end{eqnarray*}
Therefore, it follows that
\begin{eqnarray}
\int_{\bar M^*} u^2 t^{-2}\bar S_1 \dd \bar M & = & \int^1_0
t^{n-4}\sum_j b_j^2 \,\dd t \le \frac{4}{(n-3)^2}\, \int^1_0
t^{n-2}\sum_j
(b_j')^2 \dd t\nonumber\\
&  = & \frac{4}{(n-3)^2}\, \int_{\bar M^*} u_t^2 \bar S_1 \dd \bar
M. \label{c-s}
\end{eqnarray}
Finally, we conclude that
\begin{equation}
\int_{\bar M^*} u \,L u\,\dd \bar M \le
 \left(\frac{4\mu_1^-}{(n-3)^2}-1\right)
 \int_{\bar M^*} u_t^2  \bar S_1 \,\dd\bar M.
\end{equation}
Suppose $n\ge 4$ and define
\[
\mu_{\bar M}:= (1- 4\mu^-_1/ (n-3)^2).
\]
We suppose that $\mu_{\bar M}\ge 0$. Hence, it follows from
(\ref{c-s}) that
\begin{eqnarray*}
-\int_{\bar M^*} u \,L u\ge \mu_{\bar M} \, \int_{\bar M^*} u_t^2
\bar S_1 \dd \bar M \ge \mu_{\bar M} \frac{(n-3)^2}{4}\int_{\bar
M^*} u^2 t^{-2}\bar  S_1 \dd \bar M.
\end{eqnarray*}
Now, we define the truncated cone $\bar M_{\sigma,\tau}$ as the set
of points $t\theta$ in $\bar M^*$ with $0<\sigma < t < \tau\le 1$.
Let $\lambda_{\sigma,1}$ be the smallest eigenvalue of the Dirichlet
eigenvalue problem
\begin{eqnarray*}
L u +  t^{-2}\bar S_1\lambda u = 0\,\,\, \textrm{on}\,\,\, \bar
M_{\sigma,1}, \quad u = 0 \,\,\,\textrm{on}\,\,\,
\partial \bar M_{\sigma,1}.
\end{eqnarray*}
Hence, we may characterize $\lambda_{\sigma,1}$ as the Rayleigh
quotient
\begin{equation}
\lambda_{\sigma,1}= -\inf_{u\in C^1_0(\bar M_{\sigma,1}),\,
u\not\equiv 0} \frac{\int_{\bar M^*} u \, Lu\, \dd \bar M }
{\int_{\bar M^*} \frac{u^2}{t^2} \bar S_1 \,\dd \bar M}.
\end{equation}
We define
\begin{equation} I :=\inf_{u\in C^1_0(\bar M^*)} \bigg(-\int_{\bar M^*} u \, Lu \, \dd \bar M\bigg)
\end{equation}
and
\begin{equation} I_+ :=\inf_{u\in C^1_0(\bar M^*),u\not\equiv 0} \frac{-\int_{\bar M^*} u \, Lu \, \dd \bar
M}{\int_{\bar M^*} \frac{u^2}{t^2}\bar S_1 \,\dd \bar M}.
\end{equation}
Therefore, if $\mu_{\bar M}\ge 0$ (respectively, $\mu_{\bar
\mu}>0$), then $I\ge 0$ and $\inf_{\sigma} \lambda_{\sigma, 1}\ge 0$
(respectively, $I_+>0$ and $\inf_{\sigma} \lambda_{\sigma, 1}> 0$).
In the first case, we say that $\bar M^*$ is $1$-stable. In the
second case, $\bar M^*$ is said to be strictly $1$-stable.

Thus, we have proved that $\mu_{\bar M}\ge 0$ (respectively,
$\mu_{\bar M}>0$) implies that $\bar M^*$ is $1$- stable
(respectively, strictly $1$-stable).

Conversely, if  $\mu_{\bar M} <0$, then $\bar M^*$ is not
$1$-stable. In fact, in this case, we have  $\mu_1 < -(n-3)^2/4$.
Thus, the root $\gamma_1$ of $\gamma^2 + (n-3)\gamma - \mu_1 =0$ is
not real. Moreover, the function $u_1 = \Re(t^{\gamma_1}\phi_1)$ is
a Jacobi field, i.e., a solution for $\bar L_1 u - 3\bar S_3 u = 0$.
Notice that $u_1(t,\theta)=0$ for all $\theta$ whenever
$t^{\gamma_1}$ is a pure imaginary number. This happens if and only
if $\ln t \Im \gamma_1 = k\pi/2$, where $k$ is a negative integer.
Thus, we choose $\sigma, \tau$ so that
$u_1(\sigma,\cdot)=u_1(\tau,\cdot)=0$  and define the test function
for the Rayleigh quotient
\[
w(t,\theta) = u(t,\theta) \,\,\,\textrm{if}\,\,\, \sigma < t <
\tau\quad  \textrm {and} \quad w = 0 \,\,\,\textrm{otherwise}.
\]
It is clear that $w$ is a piecewise differentiable function which
satisfies
\[
\int_{\bar M^*}\big( \langle \bar T_1 \bar\nabla w, \bar\nabla w
\rangle +3\bar S_3 w \big)\dd \bar M = 0.
\]
So, $\lambda_{\sigma/2,1}<0$ since the compact support of $w$ is
strictly contained in the truncated cone $\bar
M_{\frac{\sigma}{2},1}$. We  conclude that $\inf_{\sigma}
\lambda_{\sigma,1}<0$.

In a similar way, we may prove that if $\mu_{\bar M}= 0$, then $\bar
M^*$ is not strictly $1$-stable.

These results can now be used to prove

\begin{theorem}
If $\bar M^*$ is strictly $1$-stable, then the graph $\bar
M^*_\lambda$ of the function $u_{\lambda}$ given in Theorem 1 is
strictly $1$-stable for $\lambda$ sufficiently small.
\end{theorem}

\noindent {\it Proof.} Let $\bar S_r(\lambda)$, $1\le r\le n$,
denote the elementary symmetric functions of the eigenvalues of the
Weingarten map $\bar A(\lambda)$ of $\bar M^*_\lambda$. We also
denote $\bar T_1(\lambda)=\bar S_1(\lambda)\textrm{Id}-\bar
A(\lambda)$.

As $\bar S_3(\lambda)$ depends on the Hessian of $u_\lambda$, it
follows from the $C^{2,\alpha}$ estimates on $u_\lambda$ given in
Proposition 2 that
\begin{equation}
 \sup_{\lambda}
\frac{1}{\lambda} \sup_{\bar M^*} \frac{1}{t^3} \big(\bar
S_3(\lambda) - S_3\big)<\infty.
\end{equation}
Consequently, for small $\lambda$, it holds that
\[
\int_{\bar M^*_\lambda}(\langle \bar T_1(\lambda) \bar\nabla u,
\bar\nabla u\rangle -\bar S_3(\lambda) u^2)\dd \bar M \ge \mu_{\bar
M}/2
>0,
\]
for all $u\in C^{1}_0 (\bar M^*_\lambda)$ with
$$\int_{\bar M^*_\lambda} \frac{u^2}{ t^2}\bar S_1(\lambda) \,\dd \bar M= 1.$$
This finishes the proof of the theorem.

\vspace{1.5cm}

\noindent Jorge H. S. de Lira \\
Departamento de Matem\'atica\\
Universidade Federal do Cear\'a\\
Bloco 914, Campus do Pici\\
60455-760, Fortaleza - Cear\'a, Brasil\\
jorge.lira@pq.cnpq.br

\vspace{1cm}
\noindent Marc Soret\\
Laboratoire de Math\'ematiques et Physique Th\'eorique\\
Universit\'e
de Tours\\
Parc
de Grandmont, 37200, Tours, France\\
Marc.Soret@lmpt.univ-tours.fr

\end{document}